\documentclass[12pt]{amsart}
\usepackage{amsmath,amssymb,array,longtable, MnSymbol, yfonts, amsthm}
\usepackage[mathscr]{euscript}
\usepackage{microtype}
\usepackage[colorlinks=true,urlcolor=black,citecolor=black,linkcolor=black,%
pdftitle={An analogue of Liouville's theorem and an application to cubic surfaces},%
pdfauthor={David McKinnon and Mike Roth},%
pdfsubject={Local version of Bombieri-Lang pheomena},%
pdfkeywords={Algebraic geometry, Diophantine approximation, Seshadri constants}]{hyperref}
\usepackage[colorlinks=true,urlcolor=black,citecolor=black,linkcolor=black]{hyperref}

\begin{document}
\bibliographystyle{plain}
\newtheorem{theorem}{Theorem}[section]
\newtheorem{conjecture}[theorem]{Conjecture}
\newtheorem{proposition}[theorem]{Proposition}
\newtheorem{lemma}[theorem]{Lemma}
\newtheorem{corollary}[theorem]{Corollary}
\newtheorem{definition}[theorem]{Definition}
\newtheorem{observation}[theorem]{Observation}
\def\N{{\mathbb N}}
\def\A{{\mathbb A}}
\def\F{{\mathbb F}}
\def\Q{{\mathbb Q}}
\def\R{{\mathbb R}}
\def\C{{\mathbb C}}
\def\P{{\mathbb P}}
\def\Z{{\mathbb Z}}
\def\v{{\mathbf v}}
\def\w{{\mathbf w}}
\def\x{{\mathbf x}}
\def\y{{\mathbf y}}
\def\O{{\mathcal O}}
\def\M{{\mathcal M}}
\def\E{{\mathcal E}}
\renewcommand{\AA}{\mathbb{A}}          
\newcommand{\CC}{\mathbb{C}}            
\newcommand{\NN}{\mathbb{N}}            
\newcommand{\PP}{\mathbb{P}}            
\newcommand{\QQ}{\mathbb{Q}}            
\newcommand{\RR}{\mathbb{R}}            
\renewcommand{\H}{\overline{H}}  
\def\Pos{{(Pos)}}
\newcommand{\Qbar}{\overline{k}}
\newcommand{\mult}{\operatorname{mult}}
\newcommand{\brmult}{\operatorname{brmult}}
\newcommand{\Span}{\operatorname{span}}
\newcommand{\Xtil}{\widetilde{X}}
\newcommand{\Ytil}{\widetilde{Y}}
\newcommand{\Ztil}{\widetilde{Z}}
\newcommand{\Ptil}{\widetilde{\P}}
\newcommand{\Pic}{\operatorname{Pic}}
\newcommand{\sep}{\operatorname{d_v}}
\newcommand{\sepv}{d_v}
\newcommand{\sepw}{d_w}
\newcommand{\sepwp}{d_{w'}}
\newcommand{\olddist}{\operatorname{dist}}
\newcommand{\dist}{\operatorname{dist}}
\newcommand{\squareplus}{\boxplus}
\newcommand{\PGL}{\operatorname{PGL}}
\newcommand{\Xt}{\widetilde{\X}}
\newcommand{\imXt}{\widetilde{\imX}}
\newcommand{\imX}{\textfrak{X}}     
\newcommand{\Spec}{\operatorname{Spec}}
\newcommand{\Osh}{\mathcal{O}}
\def\tr{\mbox{Tr}}
\newcommand{\bpf}{\noindent{\em Proof:} \/}
\newcommand{\epf}{\qed}
\def\qed{$\Box$}
\def\ord{\mbox{ord}}
\newcommand{\xseq}{\{x_i\}}
\newcommand{\xseqtox}{\{x_i\}\rightarrow x}
\newcommand{\zseq}{\{z_i\}}
\newcommand{\xsubseq}{\{x_i'\}}
\newcommand{\NE}{\overline{\operatorname{NE}}}
\newcommand{\st}{\,\,\rule[-0.2cm]{0.02cm}{0.6cm}\,\,}
\newcommand{\ESup}{\operatorname{ESupp}}
\newcommand{\Supp}{\operatorname{Supp}}
\newcommand{\ep}{\epsilon}
\renewcommand{\geq}{\geqslant}
\renewcommand{\leq}{\leqslant}
\newcommand{\Cnrm}[1]{||#1||}     
\newcommand{\nrm}[1]{|#1|}        
\newcommand{\uvec}{\mathbf{u}}
\newcommand{\wvec}{\mathbf{w}}
\newcommand{\fvec}{\mathbf{v}}
\newcommand{\Ish}{\mathscr{I}_{\Delta}}
\newcommand{\old}[1]{\oldstylenums{#1}}
\newcommand{\Gal}{\operatorname{Gal}} 
\newcommand{\Ctil}{\widetilde{C}}
\newcommand{\hyp}{h} 

\newcommand{\tablefont}{\tiny}
\newcounter{tablenumber}
\newcommand{\tabnum}{\refstepcounter{tablenumber}\thetablenumber}

\renewenvironment{equation}{\medskip\noindent\refstepcounter{subsection}\makebox[0pt][l]{({\bf\thesubsection})}\begin{minipage}[b]{\textwidth}$$}{$$\end{minipage}\medskip\noindent}

\renewcommand{\labelenumi}{{({\alph{enumi}})}}

\newcommand{\comment}[1]{[{\bf #1}]}

\title{An analogue of Liouville's Theorem and an application to cubic surfaces}

\author{David McKinnon}
\address{Department of Pure Mathematics \\
University of Waterloo \\
Waterloo, ON\ \  N2L 3G1 \\
Canada}
\email{dmckinnon@math.uwaterloo.ca}
\indent
\thanks{David McKinnon was partially supported by an NSERC research grant.}

\author{Mike Roth}
\address{Dept.\ of Mathematics and Statistics, Queens University, Kingston,
Ontario, Canada}
\email{mikeroth@mast.queensu.ca}
\indent
\thanks{Mike Roth was partially supported by an NSERC research grant.}

\date{Revised on January 18, 2013}

\begin{abstract}
We prove a strong analogue of Liouville's Theorem in Diophantine approximation for points on arbitrary algebraic
varieties.  We use this theorem to prove a conjecture of the first author for cubic surfaces in $\P^3$.
\end{abstract}

\subjclass[2000]{Primary}

\maketitle

\section{Introduction}

The famous theorem of K.F.\ Roth (see for example \cite[Part D]{HS}) gives a sharp upper bound on how well an
irrational algebraic number can be approximated by rational numbers.  In \cite{MR},
the authors prove an analogue of Roth's Theorem for algebraic points on arbitrary algebraic varieties.

In this paper we generalize, in the sense of \cite{MR}, Liouville's approximation theorem to arbitrary varieties,
as well as giving an extension involving the asymptotic base locus.  On $\PP^1$, except for the case that
$x\in\PP^1$ is a rational point of the base number field, Liouville's theorem is weaker than Roth's. On arbitrary
varieties the extension involving the asymptotic base locus makes it slightly more useful and we use this
to verify a conjecture of the first author for cubic surfaces in $\PP^3$.

The point of view of \cite{MR} is that the Roth and Liouville theorems are examples of
``local Bombieri-Lang phenomena'' whereby local positivity of a line bundle influences local accumulation of
rational points.  Specifically, given a variety $X$, an algebraic point $x\in X$,
and an ample line bundle $L$ on $X$, these theorems are expressed as inequalities between
$\ep_{x}(L)$, the {\em Seshadri constant}, measuring local positivity of $L$ near $x$,
and
$\alpha_{x}(L)$, an invariant measuring how well we can approximate $x$ by rational points.

\noindent
In \S\ref{sec:alpha-epsilon} we review the definitions and elementary properties of $\alpha_{x}$ and $\ep_{x}$.
In \S\ref{sec:Liouville-bound} we prove the generalized Liouville theorem (Theorem~\ref{thm:Liouville-bound}).
We close the paper in \S\ref{sec:cubic} by computing $\alpha_{x}$ and $\ep_{x}$ for an arbitrary nef line bundle
and rational point, not on a line, on a smooth cubic surface (where the lines are also rational); we then
use this to verify Conjecture~3.2 from \cite{McK}.

\section{\texorpdfstring{Elementary properties of $\alpha$ and $\epsilon$}{Elementary properties of alpha and epsilon}}\label{sec:alpha-epsilon}

In this section, we give a brief overview of the properties of $\alpha$ and $\epsilon$ used in this paper.  For a more detailed discussion of $\alpha$, see \cite{MR}.  For a more detailed discussion of $\epsilon$, there are many good references -- see for example \cite[chap.\ 5]{PAG}.  Proofs of all of the facts listed below can be found in \cite{MR}.

\noindent
{\bf The constant \boldmath{$\alpha_{x}$}.}
In order to motiviate the definition of $\alpha_{x}$ it is helpful to recall the classical case of approximation
on the line.  For a point $x\in \RR$ the {\em approximation exponent} $\tau_{x}$ of $x$ is the unique
extended real number $\tau_{x}\in(0,\infty]$ such that the inequality
$$\left|{x-\frac{a}{b}}\right| \leq \frac{1}{b^{\tau_{x}+\delta}}$$
has only finitely many solutions $a/b\in \QQ$ whenever $\delta>0$ (respectively has infinitely solutions $a/b\in \QQ$
whenever $\delta<0$).
The approximation exponent measures a certain tension between our ability to closely approximate
$x$ by rational numbers (the distance term $|x-a/b|$)  and the complexity (the $1/b$ term) of the number required to
make this approximation.  In this notation the 1844 theorem of Liouville \cite{L} is that $\tau_{x}\leq d$
for $x\in \RR$ algebraic of degree $d$ over $\QQ$.

To generalize $\tau_x$ to arbitrary projective varieties defined over a number field $k$ we replace the function $|x-a/b|$ by a distance function $\sepv(x,\cdot)$ depending on a place $v$ of $k$, and measure the complexity of a point via a height function $H_{L}(\cdot)$ depending on an ample line bundle $L$. For an introduction to the theory of heights the reader is referred to any one of \cite[Chap.\ 2]{BG}, \cite[Part B]{HS}, \cite[Chap.\ III]{La}, or \cite[Chap.\ 2]{Se}.
Unless otherwise specified all height functions in this paper are multiplicative, relative to $k$, and come from line bundles on $X$ defined over $k$. 
In this paper we normalize our height functions as follows.  The absolute values are normalized with respect to $k$: if $v$ is a finite place of $k$, $\pi$ a uniformizer of the corresponding maximal ideal, and $\kappa$ the residue field then $\nrm{\pi}_{v}=1/\#\kappa$;  if $v$ is an infinite place corresponding to an embedding $i\colon k\hookrightarrow \CC$ then $\nrm{x}_{v} = \nrm{i(x)}^{m_v}$ for all $x\in k$, where $m_v=1$ or $2$ depending on whether $v$ is real or complex.  The heights are then normalized so that for a point ${x}=[x_0:\cdots :x_n]\in\PP^n(k)$, the height with respect to $\Osh_{\PP^n}(1)$ is \[H({x})=\prod_v \max(\nrm{x_0}_{v},\ldots, \nrm{x_n}_{v})\] where the product ranges over all the places $v$ of $k$.

In order to define a distance function we fix a place $v$ of $k$ and extension (which we also call $v$) to $\Qbar$.

\noindent
{\bf If $v$ is archimedean:}
We choose a distance function
on $X(\Qbar)$ by choosing an embedding $X\hookrightarrow\PP^{n}_{k}$ defined over $k$, and
pulling back (via $v$) the distance function on $\PP^{n}(\CC)$ given by the Fubini-Study metric on $\PP^{n}$.
We denote this distance by $\dist(\cdot,\cdot)$.   We set $\sepv(\cdot,\cdot)=\dist(\cdot,\cdot)^{m_v}$
where $m_v=1$ if $v$ is real and $m_v=2$ if $v$ is complex.
This distance function depends on the choice
of embedding, but by \cite[Proposition 2.1]{MR} any two embeddings give equivalent distance functions and the
choice of embedding will not matter for the definition of $\alpha_{x}$.

\noindent
{\bf If $v$ is non-archimedean:}
Again choose a projective embedding $X\hookrightarrow \PP^{m}$ defined over $k$.
If $x,y\in X(\Qbar)$, consider the corresponding projective
coordinates $\mathbf{x}=[x_0\colon\cdots\colon x_m]$, $\mathbf{y}=[y_0\colon\cdots\colon y_m]$, and set
$\sepv(x,y) = \frac{H_v(\mathbf{x}\wedge\mathbf{y})}{H_{v}(\mathbf{x})H_{v}(\mathbf{y})}$ where $H_v$ is the local
height at the place $v$ (this is the definition given in
\cite[2.8.16]{BG} although we are using a different normalization for height than \cite{BG}).

This definition is somewhat opaque on first reading but is a compact way of stating a very concrete notion of
$v$-adic distance:
points $x$ and $y$ are close if the corresponding curves in an integral model of $X$ have high order of contact at
the place $v$ (see e.g., \cite[\S2]{MR}).  In other words, two points $x$ and $y$ are close if they are congruent modulo a high power of the maximal ideal $\mathfrak{m}_v$ of $\O_{k_v}$.  For any fixed $x\in X(\Qbar)$, different embeddings give equivalent
functions $\sepv(x,\cdot)$, see \cite[Corollary 2.3]{MR}.

\begin{definition}\label{def:alpha}
Let $X$ be a projective variety defined over a number field $k$, $L$ an ample line bundle defined over $k$,
and $x\in X(\Qbar)$.  Then we define $\alpha_x=\alpha_{x}(L)$ to be the unique extended real number 
$\alpha_x\in(0,\infty]$ such that the inequality
$$\sepv(x,y)^{\alpha_{x}+\delta}< H_{L}(y)^{-1}$$
has only finitely many solutions $y\in X(k)$ (respectively infinitely many solutions $y\in X(k)$) for any $\delta<0$
(respectively any $\delta>0$).
\end{definition}

The one essential change in our definition of $\alpha_{x}$ over $\tau_{x}$ is that we have moved the exponent from
the height term to the distance term.  As a result, for $x\in\RR=\AA^1(\RR)\subset\PP^1(\RR)$ we have
$\alpha_{x}(\Osh_{\PP^1}(1))=\frac{1}{\tau_{x}}$.  In particular
the theorem of Liouville becomes $\alpha_{x}(\Osh_{\PP^1}(1))\geq \frac{1}{d}$
for $x\in \RR$ of degree $d$ over $\QQ$, and it is this type of lower bound
that we wish to generalize to arbitrary varieties.
The choice of moving the exponent is justified by the resulting formal similarity with the Seshadri constant,
and more natural behaviour when we vary $L$ (see, for example, Proposition \ref{prop:alpha-and-ep}).

In proving results about $\alpha_{x}$ it is useful to have a characterization of $\alpha_{x}$ in terms of
``test sequences'', and to associate an approximation constant to such a sequence.

\begin{definition}\label{seqappconst}
Let $X$ be a projective variety, $x\in X(\Qbar)$, $L$ a line bundle
on $X$.  For any sequence $\xseq\subset X(k)$ of distinct
points with $\sepv(x,x_i)\rightarrow 0$ (which we denote by $\xseqtox$), we set

$$A(\xseq, L) = \left\{{
\gamma\in\R \st
\sepv(x,x_i)^{\gamma} H_{L}(x_i)\,\,\mbox{is bounded from above}
}\right\}.
$$
\end{definition}

\noindent
{\bf Remarks.}
(a) It follows easily from the definition that if $A(\xseq, L)$ is nonempty then it is an
interval unbounded to the right, i.e., if $\gamma\in A(\xseq,L)$ then $\gamma+\delta\in A(\xseq,L)$ for any $\delta>0$.

(b) If $\xsubseq$ is a subsequence of $\xseq$ then $A(\xseq,L)\subseteq A(\xsubseq,L)$.

\begin{definition}
If $A(\xseq,L)$ is empty we set $\alpha_{x}(\xseq,L)=\infty$.  Otherwise we set
$\alpha_{x}(\xseq, L)$ to be the infimum of $A(\xseq,L)$.  We call $\alpha_{x}(\xseq,L)$ the approximation constant
of $\xseq$ with respect to $L$.
\end{definition}

As $i\to\infty$ we have $\sepv(x,x_i)\to0$ and $H_{L}(x_i)\to\infty$.  We thus expect that
$\sepv(x,x_i)^{\gamma}H_{L}(x_i)$ goes to $0$ for large $\gamma$ and to $\infty$ for small $\gamma$.  The number
$\alpha_{x}(\xseq,L)$ marks the transition point between these two behaviours.

By remark (b) above if $\xsubseq$ is a subsequence of $\xseq$ then $\alpha_{x}(\xsubseq,L)\leq \alpha_{x}(\xseq,L)$.
Thus we may freely replace a sequence with a subsequence when trying to establish lower bounds.

\begin{proposition}
Let $X$ be a projective variety defined over a number field $k$, $L$ an ample line bundle defined over $k$,
and $x\in X(\Qbar)$.  Then $\alpha_{x}(L)$ is the infimum of of all approximation constants of
sequences of points in $X(k)$ converging to $x$.  If no such sequence exists then $\alpha_{x,X}(L)=\infty$.
\end{proposition}

\bpf
This is an elementary argument using sequences and the fact that if $L$ is ample there are only finitely many
rational points of bounded height.  For details see \cite[Proposition 2.9]{MR}.
\epf

The following lemma
gives an equivalent local expression for the distance, useful for calculating with test sequences.

\begin{lemma}
\label{lem:local-distance}
Let $x$ be a point of $X(k)$ and let $U$ be an open affine of
$X$ containing $x$.
Let $u_1$, \ldots, $u_r$ be elements of $\Gamma(U,\Osh_{X})$
which generate the maximal ideal of $x$.  Then
there are constants $c$ and $C$ such that
$$c\,\sepv(x,y)\leq \min\left({1,\max\left(\nrm{u_1(y)}_v,\ldots, \nrm{u_r(y)}_v\right)\rule{0cm}{0.4cm}}\right)
\leq C\, \sepv(x,y)$$
for all $y\in U(\Qbar)$.  I.e., on $U(\Qbar)$ the function
$\min(1,\max\left(\nrm{u_1(\cdot)}_v,\ldots, \nrm{u_r(\cdot)}_v\right))$
is equivalent to the function $\sepv(x,\cdot)$.
\end{lemma}

\bpf See \cite[Lemma 2.4]{MR}. \epf

We need two results on $\alpha_{x}$ before continuing onto the Seshadri constant.
First, we will need to know how to calculate $\alpha_{x}$ in one simple case.

\begin{lemma}\label{lem:projective}\label{prop:Roth-on-P1}
Let $x\in\P^n(k)$.  Then $\alpha_{x,\P^n}(\O_{\P^n}(1))=1$.
\end{lemma}

\bpf
This is Lemma~{2.11} from \cite{MR}.  \epf

Second, it will be useful to know how the approximation constant changes when we change the field $k$.
We use the notation that for an extension field $K/k$, $\alpha_{x}(\xseq,L)_{K}$
(respectively $\alpha_{x}(L)_{K}$) denotes the
approximation constant of a sequence (resp.\ point $x$) computed with respect to $K$.   This means that when
computing $\alpha$, we use the height $H_{L}$ relative to $K$ and normalize $\sepv$ relative to $K$.
If $d=[K\colon k]$ and $m_v=[K_v\colon k_v]$
(where $K_v$ and $k_v$ denote the completions of $K$ and $k$ with respect to $v$)
then this means simply that $H_{L}(x_i)_{K} = H_{L}(x_i)^{d}_{k}$
and $\sepv(x,x_i)_{K} = \sepv(x,x_i)^{m_v}_{k}$.

\begin{proposition}\label{prop:change-of-field}
Suppose $x\in X(\Qbar)$, $L$ a line bundle defined over $k$, and $\xseqtox$ a sequence of points in $X(k)$
approximating $x$.   Let $K$ be any finite extension of $k$.  Then $\xseqtox$ can also be considered
to be a set of points of $X(K)$ approximating $x$.   Set $m_v=[K_v\colon k_v]$, and let $d=[K\colon k]$.
Then
$$\alpha_x\left({\xseq, L}\right)_{K} =
\frac{d}{m_v} \alpha_x\left({\xseq, L}\right)_{k}.$$
In particular, we have the bound $\alpha_x(L)_{K} \leq \frac{d}{m_v} \alpha_x(L)_{k}$.
\end{proposition}

\bpf
The claim that
$\alpha_x\left({\xseq, L}\right)_{K} =
\frac{d}{m_v} \alpha_x\left({\xseq, L}\right)_{k}$ follows immediately from the equalities
$H_{L}(\cdot)_{K} = H_{L}(\cdot)_{k}^{d}$
and $\sepv(\cdot,\cdot)_K=\sepv(\cdot,\cdot)_k^{m_{v}}$.  The inequality
$\alpha_x(L)_{K} \leq \frac{d}{m_v} \alpha_x(L)_{k}$ then follows since the sequences of $k$-points 
approximating $x$ are a subset of the sequences of $K$-points approximating $x$.
\epf

\noindent
{\bf Remark.} \label{rem:same-local-degree}
Let $x$ be a point of $X(\Qbar)$ and let $K$ be the field of definition of $x$.
If $K\not\subseteq k_v$, or equivalently, $K_v\neq k_v$ then it will be impossible to find a sequence of
points of $X(k)$ converging (in terms of $\sepv$) to $x$. For example, when $v$ is archimedean this happens when
$k_v=\RR$ and $K_v=\CC$.
Thus, if we can approximate $x$ by points of $X(k)$ we may assume that $K_v=k_v$ and so $m_v=1$.

The following result (appearing in \cite{MR} as Theorem 2.14, and incorrectly in \cite{McK} as Theorem~2.8) is
obtained by combining the Roth and Dirichlet theorems for approximation on $\PP^1$,
as well as the local information about the singularity type, shows how to calculate $\alpha_{x}$ on any
singular $k$-rational curve.

\begin{theorem}\label{thm:curve}
Let $C$ be any singular $k$-rational curve and $\varphi\colon\PP^1\rightarrow C$ the normalization map.
Then for any ample line bundle $L$ on $C$, and any $x\in C(\Qbar)$ we have the equality:
\[\alpha_{x,C}(L)=\min_{q\in \varphi^{-1}(x)} d/r_{q} m_{q}\]
where $d=\deg(L)$, $m_{q}$ is the multiplicity of the branch of $C$ through $x$ corresponding to $q$, and
\[r_{q}=
\begin{cases}
0 & \text{if $\kappa(q)\not\subseteq k_v$} \\
1 &\text{if $\kappa(q)=k$} \\
2 &\text{otherwise.}
\end{cases}
\]
\end{theorem}

\medskip
\noindent
Here $\kappa(q)$ means the residue field of the point $q$, and we use $r_q=0$ as a shorthand for $d/r_{q}m_{q}=\infty$.

\noindent
{\bf The Seshadri constant.}
The Seshadri constant was introduced by Demailly in \cite{Dem} for the purposes of measuring the local positivity
of a line bundle.

\begin{definition}\label{def:seshadri}
Let $X$ be a  projective variety, $x$ a point of $X$, and $L$ a
nef line bundle on $X$.  The {\em Seshadri constant},
$\ep_{x}(L)$, is defined to be

\[\ep_{x}(L) :=
\sup\left\{{\gamma\geq 0 \mid \pi^{*}L -\gamma E\,\,\,
\mbox{is nef}\,}\right\}\]
where $\pi:\Xtil\longrightarrow X$ is the blowup of $X$ at $x$, with
exceptional divisor $E$.
\end{definition}

In the discussion of Conjecture \ref{conj:ratcurve} below we will need the following alternate
characterization of the Seshadri constant:

\begin{proposition}\label{prop:seshcurvemult}
With the same setup as definition \ref{def:seshadri},
\[
\ep_{x}(L) =
\inf_{x\in C\subseteq\, X}
\left\{
\frac{(L\cdot C)}{\mult_x(C)}
\right\}
\]
where the infimum is taken over all reduced irreducible curves $C$ passing
through $x$.
\end{proposition}

\bpf
This is \cite[Proposition 5.15]{PAG}.  \epf

In order to indicate the parallels between $\alpha_{x}$ and $\ep_{x}$, and for use below,
we list a few of their formal properties here.

\begin{proposition}\label{prop:alpha-and-ep}
Let $X$ be a projective variety defined over $k$, $x\in X(\Qbar)$, and let $L$ be any ample line bundle on $X$
(also defined over $k$, following our conventions above).

\begin{enumerate}
\item For any positive integer $m$, $\alpha_{x}(mL) = m\alpha_{x}(L)$ and $\ep_{x}(mL)=m\ep_{x}(L)$.
(Thus $\alpha$ and $\ep$ also make sense for ample $\QQ$-divisors.)
\item  $\alpha_{x}$ and $\ep_{x}$ are concave functions of $L$:
for any positive rational numbers $a$ and $b$, and any ample $\Q$-divisors $L_1$ and $L_2$  (again defined over $k$)
we have
$$\alpha_x(a L_1+ b L_2) \geq a\alpha_x(L_1)+b\alpha_x(L_2)
\,\,\,\mbox{and}\,\,\,
\ep_x(a L_1+ b L_2) \geq a\ep_x(L_1)+b\ep_x(L_2)$$
\item If $Z$ is a subvariety of $X$ defined over $k$ then for any point $z\in Z(\Qbar)$
we have $\alpha_{z}(L|_Z)\geq \alpha_{z,X}(L)$ and $\ep_{z}(L|_{Z})\geq \ep_{z}(L)$.
\item If $Y$ is also a variety defined over $k$, $x\in X(k)$, $y\in Y(k)$ and $L_X$ and $L_Y$ are nef
line bundles defined on $X$ and $Y$ respectively then
\[\alpha_{x\times y,X\times Y}(L_X\squareplus L_Y) = \min(\alpha_{x,X}(L_X), \alpha_{y,Y}(L_Y))\]
and
\[\ep_{x\times y,X\times Y}(L_X\squareplus L_Y) =
\min(\ep_{x,X}(L_X), \ep_{y,Y}(L_Y)).\]
\end{enumerate}
\end{proposition}

Note that by $L_{X}\squareplus L_Y$ we mean the line bundle
$pr_{X}^{*}L_1 + pr_{Y}^{*}L_2$ on $X\times Y$, where $pr_{X}$ and $pr_{Y}$ are the projections.
We prefer additive notation for line bundles since this is in line with the behaviour of $\alpha_{x}$ and $\ep_{x}$,
and hence use $L_{X}\squareplus L_Y$ rather than $L_1\boxtimes L_2$.

\bpf
All the proofs follow from elementary arguments using the definitions.
For the statements about $\alpha_{x}$ see
\cite[Proposition 2.12]{MR}, and for the statements about $\ep_{x}$ see \cite[Proposition 3.4]{MR}.  
\epf

\section{\texorpdfstring{A Liouville lower bound for $\alpha$}{A Liouville lower bound for alpha}}\label{sec:Liouville-bound}

In this section, as in the previous one, we fix a number field $k$ and let $X$ be a projective variety defined over $k$.

%

\begin{lemma}\label{lem:height-of-E}
Let $x$ be a point of $X(k)$, and $\pi\colon\Xtil\longrightarrow X$ the blow up of $X$ at $x$ with exceptional
divisor $E$.
Choose an embedding $\varphi\colon X\hookrightarrow \P^n$ so that $x\mapsto [1\colon 0\colon\cdots\colon 0]$.
Let $Z_0$,\ldots, $Z_n$ be the coordinates on $\P^n$ and define functions $u_i$, $i=1,\ldots, n$ on the
open subset where $Z_0\neq 0$ by $u_i=Z_i/Z_0$.

For each place $w$ of $k$, define a function $e_w\colon X(k)\rightarrow \R_{\geq 0}$ by

$$e_w(y) = \left\{{\begin{array}{cl}
1 & \mbox{if $Z_0(y)=0$,} \\
\min\left({1, \max(\nrm{u_1(y)}_{w},\ldots, \nrm{u_n(y)}_{w})}\right) & \mbox{if $Z_0(y)\neq 0$.} \\
\end{array}}\right.
$$

Then

\begin{enumerate}
\item $e_w\leq 1$ for all places $w$.
\item $e_v$ is equivalent to $\sepv$.
\item For $y\in X(\Qbar)$, $y\neq x$,  we have
$H_{E}(y) = \left(\prod_{w} e_{w}(y)\right)^{-1}$.
\end{enumerate}
\end{lemma}

\bpf
Part (a) is clear from the definition.  Part (b) is precisely
Lemma \ref{lem:local-distance}.
In (c) we are considering points $y\in X(k)$, $y\neq x$ also to be points of $\Xtil(k)$ via the birational
map $\pi$.
To prove (c) it suffices, by using the functoriality of heights under pullback,  to consider the case that $X=\P^n$.
Then the blow up $\Ptil^n$ of $\P^n$ at $x$ is a
subvariety of $\P^n\times \P^{n-1}$ and $\O_{\Ptil^n}(E)$ is the restriction of $\O_{\P^n\times\P^{n-1}}(1,-1)$
to $\Ptil^n$.  From this description of $\O_{\Ptil^n}(E)$ we obtain the formula

$$H_{E}(y) = \prod_{w}
\frac{\max\left({\nrm{Z_0(y)}_w, \nrm{Z_1(y)}_w,\ldots, \nrm{Z_n(y)}_w}\right)}%
{\max\left({\nrm{Z_1(y)}_w,\ldots, \nrm{Z_n(y)}_w}\right)}
$$
from which (c) follows easily.  \epf

\begin{lemma}\label{lem:gamma-bound}
Suppose that $x\in X(k)$ and let $\pi\colon\Xtil\longrightarrow X$ be the blow up at $x$ with exceptional divisor $E$.
Let $L$ be an ample line bundle on $X$ and $\gamma>0$ a rational number such that $L_{\gamma}:=\pi^{*}L-\gamma E$
is in the effective cone of $\Xtil$.
Let $B'$ be the asymptotic base locus of $L_{\gamma}$ and set $B=\pi(B')$.

\noindent Then for
any sequence $\xseqtox$ such that all points of $\xseq$ are outside of $B$, $\alpha(\xseq,L)\geq \gamma$.
\end{lemma}

\bpf
Let $U=\Xtil\setminus B'$.  Since $B'$ is the asymptotic base-locus of $L_{\gamma}$ there is a constant $c$
so that $H_{L_{\gamma}}(y)\geq c$ for all $y\in U(k)$.   Applying Lemma \ref{lem:height-of-E} we then have
$$c\leq H_{L_{\gamma}}(x_i) = H_{L}(x_i) H_{E}(x_i)^{-\gamma}
\stackrel{\ref{lem:height-of-E}(c)}{=}  H_{L}(x_i)\left(\prod_{w}e_{w}(x_i)\right)^{\gamma}
\stackrel{\ref{lem:height-of-E}(a)}{\leq} H_{L}(x_i) e_{v}(x_i)^{\gamma}.$$
By Lemma \ref{lem:height-of-E}(b) $\sepv(x,x_i)$ and $e_{v}(x_i)$ are equivalent functions on $X(k)$ and therefore
$H_{L}(x_i)\sepv(x,x_i)^{\gamma}\geq c'$ for some positive constant $c'$.

For any $\delta>0$ we thus have $H_{L}(x_i)\sepv(x,x_i)^{\gamma-\delta}\geq c'\sepv(x,x_i)^{-\delta}$ and
so conclude that $\gamma-\delta\not\in A(\xseq,L)$ since $c'\sepv(x,x_i)^{-\delta}\to\infty$ as $i\to\infty$. Therefore
$\gamma\leq \alpha(\xseq,L)$.
\epf

\vspace{.1in}

The main result of this section is the following implication of Lemma \ref{lem:gamma-bound}.

\begin{theorem}\label{thm:Liouville-bound}
Let $X$ be an algebraic variety defined over $k$,
$x\in X(\Qbar)$ any point, and set $d=[K:k]$ where $K$ is the field of definition of $x$.

Let $\pi\colon \tilde{X}\to X$ be the blowup of $X$ at $x$, with
exceptional divisor $E$,  $L$ an ample line bundle on $X$, and
$\gamma>0$ a rational number such that $L_{\gamma}:=\pi^{*}L - \gamma E$ is in the effective cone
of $\Xtil$.  Finally let $B'$ be the asymptotic base locus of
$L_{\gamma}$ and set $B=\pi(B')$.  Then

\begin{enumerate}
\item For any sequence $\xseqtox$ of $k$-points approximating $x$ if infinitely many points of $\xseq$ are
outside $B$ then $\alpha(\xseq,L) \geq \gamma/d$.

\medskip

\item If $\alpha_{x}(L) < \gamma/d$ then $x\in B$ and $\alpha_{x}(L) = \alpha_{x}(L|_B)$.

\medskip

\item If $x\in B$ and $\alpha_{x}(L|_{B})\geq \gamma/d$ then $\alpha_{x}(L)\geq \gamma/d$.
\end{enumerate}
\end{theorem}

Note that $\pi$, $\Xtil$, $E$, and $B'$ are only defined over $K$.  However since $\pi$ is a morphism of $k$-schemes,
$B$ is defined over $k$.

\bpf
Let $\{x_i\}$ be a sequence  approximating $x$.  If infinitely many $x_i$ lie outside of $B$ then we may pass
to the subsequence of points outside of $B$, which could only have the effect of lowering the approximation
constant of the sequence.  To prove part (a) we may therefore assume that all points of $\xseq$ lie outside $B$.
Applying Lemma~\ref{lem:gamma-bound} to estimate the approximation constant computed relative to $K$
we conclude that $\alpha(\xseq,L)_{K} \geq \gamma$.  Since there is a sequence of $k$-points approximating $x$
we conclude by the remark on page \pageref{rem:same-local-degree} that
(in the notation of Proposition~\ref{prop:change-of-field}) $m_v=1$.
Therefore by Proposition~\ref{prop:change-of-field} $\alpha(\xseq,L)_{k} = \frac{1}{d}\alpha(\xseq,L)_{K}\geq \gamma/d$,
proving (a).

If $\alpha_{x}(L)< \gamma/d$ then there must be a sequence $\xseq$ approximating $x$ such that
$\alpha(\xseq,L)< \gamma/d$.  By part (a) this implies that all but finitely many $x_i$ lie in $B$.  Thus
$x\in B$ since $B$ is closed.   Since omitting finitely many elements of a sequence does not change the approximation
constant we may assume that all $x_i$ are contained in $B$.  Since $\alpha_{x,X}(L)$ is the infimum of the approximation
constants for sequences $\xseq$ with $\alpha(\xseq,L)< \gamma/d$ we conclude that
$\alpha_{x}(L) = \alpha_{x}(L|_{B})$ proving (b).

If $\alpha_{x}(L) < \gamma/d$ then part (b) along with the hypothesis for part (c) lead to an immediate contradiction.
Thus, under the hypotheses of part (c), $\alpha_{x}(L)\geq \gamma/d$. \epf

\noindent
{\bf Remark.} Theorem~\ref{thm:Liouville-bound} still holds if we replace $B$ by the Zariski closure of $B(k)$.
This has the added advantage that every component of $B$ is then absolutely irreducible (see \cite[Lemma 2.15]{MR}).

\begin{corollary}\label{cor:alpha-vs-epsilon}
For all ample line bundles $L$ on $X$ we have $\alpha_{x}(L) \geq \ep_{x}(L)/d$.
\end{corollary}

\bpf
Let $\pi\colon\Xtil\longrightarrow X$ be the blow up of $X$ at $x$.  By the definition of $\ep_{x}(L)$ for all
rational $\gamma$ satisfying $0<\gamma< \ep_{x}(L)$ the line bundle $\pi^{*}L-\gamma E$ is ample on $\Xtil$ and
in particular the asymptotic base locus of $\pi^{*}L-\gamma E$ is empty.  Thus by
Theorem~\ref{thm:Liouville-bound}(a) we conclude that $\alpha_{x}(L)\geq\gamma/d$ for any such $\gamma$, and hence
that $\alpha_{x}(L)\geq \ep_{x}(L)/d$.
\epf

\vspace{.1in}

\noindent
{\bf Remark.} If $X=\P^1$ then Corollary~\ref{cor:alpha-vs-epsilon}
and the fact that $\ep_{x}(\O_{\P^n}(1))=1$ give $\alpha_{x}(\O_{\P^1}(1))\geq 1/d$.
Thus on $\PP^1$ Corollary~\ref{cor:alpha-vs-epsilon} amounts to the classic Liouville bound $\tau_x \leq d$.
For this reason we consider Theorem~\ref{thm:Liouville-bound} and Corollary~\ref{cor:alpha-vs-epsilon} to be
``Liouville bounds'' for $\alpha_x$.

The effective cone is usually larger than the ample cone, and
in general the parts of Theorem~\ref{thm:Liouville-bound}
imply a much stronger lower bound for $\alpha_{x}(L)$ than Corollary~\ref{cor:alpha-vs-epsilon}.
We will use this in the next section to compute $\alpha$ for the cubic surface, but give a brief illustration
now by calculating $\alpha$ for rational points of a non-split quadric surface in $\P^3$.
(For a split quadric surface $\alpha_{x}(\Osh_{\PP^1\times\PP^1}(a,b))=\min(a,b)$ when $a,b>0$, as 
implied by Proposition \ref{prop:alpha-and-ep}(d) and computed
in both \cite[Theorem 3.1]{McK} and \cite[\S2; Example (c)]{MR}.)

\vspace{.1in}

\noindent
{\bf Example.}
Let $X$ be a smooth quadric surface in $\P^3$ defined over $k$, and set $L=\O_{\P^3}(1)|_{X}$.   We
assume that no lines on $X$ are defined over $k$.
Let $x$ be a $k$-point of $X$.
By intersecting with a (rationally defined) hyperplane we may find a conic $C$ passing through $x$ such that
$C$ is isomorphic to $\P^1$ over $k$.  By Lemma~\ref{lem:projective} and
Proposition~\ref{prop:alpha-and-ep}(a,c), we therefore have
$\alpha_{x,X}(L) \leq \alpha_{x,C}(L|_{C}) = \alpha_{x,\P^1}(\O_{\P^1}(2)) = 2$.
Since $x$ lies on a line (over $\Qbar$), we have $\ep_{x}(L)=1$, and applying
Corollary~\ref{cor:alpha-vs-epsilon}
we obtain $\alpha_{x}(L)\geq 1$.  Thus $1\leq \alpha_{x}(L) \leq 2$, i.e., Corollary~\ref{cor:alpha-vs-epsilon}
does not give enough information to determine $\alpha_{x}(L)$ in this case.

However, let $\pi\colon\Xtil\longrightarrow X$ the blow up of $X$ at $x$ with exceptional
divisor $E$.  Then $\pi^{*}L-2E$ is effective with base locus
the proper transform of the two lines passing through $x$.  In particular the image $B$ of this base locus
is the union of the two lines of ruling passing through $x$.
Since (by assumption) neither of these lines is defined over $k$, $x$ is the only $k$-point of $B$.  Thus by
Theorem~\ref{thm:Liouville-bound}(a) if $\xseq$ is any sequence of $k$-points approximating $x$ then
$\alpha(\xseq,L)\geq 2$, and in particular $\alpha_{x}(L)\geq 2$.
Thus $\alpha_{x}(L)=2$ for all $k$-points of $X$.

Since $X$ is non-split the Picard group of $X$ (over $k$) has rank one with generator $L$.  Thus
the above computation and the homogeneity in Proposition \ref{prop:alpha-and-ep}(a)
determines $\alpha$ for all $x\in X(k)$ and all ample line bundles on $X$ defined over $k$.

\section{The cubic surface}\label{sec:cubic}

In this section, we will compute $\alpha_x$ and $\ep_x$ for all $k$-rational points $x$ on the
blowup $X$ of $\P^2$ at six $k$-rational points in general position.

To begin, we will recall some notions from \cite{McK}.

\begin{definition}\label{bestappcurve}
A sequence $\xseqtox$ whose approximation constant is equal to $\alpha_x(L)$
(if such a sequence exists) is called a sequence of best approximation to $x$.
A curve $C$ passing through $x$ is a called a curve of best approximation (with respect to $L$) if
$C$ contains a sequence of best approximation to $x$.
\end{definition}

\noindent
In other words, if $C$ is a curve of best
approximation to $x$ on $X$, then the rational points on $C$
approximate $x$ roughly as well as the rational points on $X$
approximate $x$.

In the example of the non-split quadric --- and in many others considered in \cite{McK} ---
there is always a curve of best approximation to $x$.
In \cite[\S4]{McK}  it is shown that if Vojta's main conjectures are true,
then $\alpha_{x}(L)$ finite implies that  $\alpha_{x}(L)$  is computed on a subvariety $V\subseteq X$
of negative Kodaira dimension (possibly $X$ itself, if $X$ has negative Kodaira dimension).
Since varieties of negative Kodaira dimension are (again, conjecturally) covered by rational curves,
one is led to the following further prediction (\cite[Conjecture 2.7]{McK}):

\begin{conjecture}\label{conj:ratcurve}
Let $X$ be an algebraic variety defined over $k$, and $L$ any ample
divisor on $X$.  Let $x$ be any $k$-rational point on $X$ and assume
that there is a rational curve defined over $k$ passing through $x$. Then there exists a curve $C$ (necessarily rational) of
best approximation to $x$ on $X$ with respect to $L$.
\end{conjecture}

In \cite{McK}, the first author proves this conjecture in many cases, and shows that in many others it follows from Vojta's Conjecture.
Those proofs use a slightly different definition of $\alpha$, but the proofs do not essentially
change in the new setting.


The Seshadri-constant analogue of a curve of best approximation is called a {\em Seshadri curve} (cf.\
Proposition~\ref{prop:seshcurvemult}):

\begin{definition}
Let $L$ be a nef divisor on an algebraic variety $X$, and $x\in X$ any
point.  A Seshadri curve for $x$ with respect to $L$ is a curve $C$ such
that $\ep_{x,X}(L)=(L\cdot C)/\mult_x(C)$.
\end{definition}

In all currently known examples, there exists a Seshadri curve for
$x$ with respect to $L$, but it is conjectured that this is not always
the case.  In particular, it is possible that the Seshadri constant
might sometimes be irrational (see \cite[Remark 5.1.13]{PAG}).

It is useful to know that for a fixed curve $C$, the set of line bundles for which $C$ is a curve
of best approximation form a subcone of the N\'eron-Severi group, and similarly for the property of
being a Seshadri curve.

\begin{proposition}\label{prop:alphacurve}
Let $X$ be a variety defined over $k$, and let $x\in X(k)$ be any $k$-rational point.  Let $D_1$ and $D_2$ be nef
divisors on $X$ with height functions $H_1$ and $H_2$ bounded below by a positive constant in some neighbourhood
of $x$.  Let $a_1$ and $a_2$ be non-negative integers, and let $D=a_1D_1+a_2D_2$.

\begin{enumerate}
\item If $C$ is a curve of best approximation for $D_1$ and $D_2$, then $C$ is also a curve of best approximation
for $D$.

\item If $C$ is a Seshadri curve for $x$ with respect to $D_1$ and $D_2$, then $C$ is also a Seshadri curve for $x$
with respect to $D$.
\end{enumerate}
\end{proposition}

\bpf
Part (a) appears as \cite[Corollary 3.2]{McK}.  To prove part (b), note that Proposition \ref{prop:alpha-and-ep}(b) implies the estimate
$$\ep_{x}(a_1 D_1 + a_2 D_2) \geq a_1 \ep_{x}(D_1) + a_2 \ep_{x}(D_2).$$
On the other hand, the hypotheses of part (b) give
$$\frac{C\cdot D}{\mult_{x}C} =
\frac{C\cdot(a_1D_1+a_2D_2)}{\mult_{x}C} =
\frac{a_1(C\cdot D_1)}{\mult_{x}C} +
\frac{a_2(C\cdot D_2)}{\mult_{x}C} =
a_1\ep_{x}(D_1)+a_2\ep_{x}(D_2).
$$
Thus, by  Proposition \ref{prop:seshcurvemult},
$a_1\ep_{x}(D_1)+a_2\ep_{x}(D_2)$ is an upper bound for $\ep_{x}(D)$.  Therefore
$\ep_{x}(a_1 D_1 + a_2 D_2) = a_1 \ep_{x}(D_1) + a_2 \ep_{x}(D_2)$ and $C$ is a Seshadri curve for $D$, proving (b).
\epf

%

We are now ready to begin the proof of the main result of this section.
Before we state and prove the general result, we will illustrate the
fundamental techniques in the case $L=-K$.

\begin{theorem}\label{thm:cubic}
Let $X$ be a smooth cubic surface in $\P^3$ defined over $k$, and
isomorphic over $k$ to the blowup of $\P^2$ at six $k$-rational points
in general position.  Let $x\in X(k)$ be any $k$-rational point, and
let $C_x$ be the curve of intersection of $X$ with the tangent plane
to $X$ at $x$.    Then
$$
\ep_{x}(-K) =
\begin{cases}
1 & \mbox{if $x$ lies on one of the $27$ lines of $X$} \\
\frac{3}{2} & \mbox{otherwise}\\
\end{cases}
$$
while
$$
\alpha_{x}(-K)  =
\begin{cases}
1 & \mbox{if $x$ lies on one of the $27$ lines of $X$} \\
\\
\frac{3}{2} & \mbox{ if $x$ is not on one of the $27$ lines, and if either}   \\
&
\mbox{\begin{minipage}{0.6\textwidth}
\noindent
\begin{itemize}
\item[$\circ$] $C_x$ is cuspidal at $x$, or
\item[$\circ$] $C_x$ is nodal at $x$ with tangent lines having slopes in $k_v$ but not $k$
\end{itemize}
\end{minipage}}
\\
\\
2 & \mbox{otherwise} \\
& \mbox{\begin{minipage}{0.6\textwidth}
(i.e., 
$C_x$ is nodal at $x$, and the slopes of the
tangent lines are in $k$ or not in $k_v$.)
\end{minipage}} \\
\end{cases}
.
$$

\end{theorem}

\noindent
{\it Proof:} Set $L=-K=\Osh_{\PP^3}(1)|_{X}$, and let $x$ be a point of $X(k)$.
If $x$ lies on a line $\ell$ then by
Proposition \ref{prop:alpha-and-ep}(c) we have
$\ep_{x,\ell}(L|_{\ell}) \geq \ep_{x,X}(L) \geq \ep_{x,\PP^3}(\Osh_{\PP^3}(1)).$
Since $\ep_{x,\ell}(L|_{\ell})=\ep_{x,\PP^3}(\Osh_{\PP^3}(1))=1$,
we conclude that $\ep_{x}(L)=1$.  Similarly (again using Proposition \ref{prop:alpha-and-ep}(c)) we conclude
that $\alpha_{x}(L)=1$.

We now suppose that $x$ does not lie on a line.
Let $\pi\colon Y\to X$ be the blowup of $X$ at $x$, with exceptional divisor $E$.  Then $C_x$ is a
Seshadri curve for $x$ with respect to $L$.  To see this, note first
that $C_x$ satisfies $C_x.L/\mbox{mult}_x(C_x)=3/2$, so
$\ep_x(L)\leq 3/2$.  Conversely, if $a>3/2$, then $\pi^*L-aE$ is
not nef, because $(\pi*L-aE)(\pi^*L-2E) =3-2a<0$ and $\pi^*L-2E$ is
the class of the proper transform of $C_x$.  Thus, $\ep_x(L)\geq 3/2$,
implying $\ep_x(L)=3/2$, and $C_x$ is a Seshadri curve for $x$
with respect to $L$.

We now turn to the computation of $\alpha$.
The asymptotic base locus of $\pi^{*}L-2E$ is $\Ctil_{x}$, the proper transform of $C_x$.
Hence by Theorem \ref{thm:Liouville-bound}(b) either $\alpha_{x}(L)\geq 2$ or
$\alpha_{x}(L) = \alpha_{x,C_x}(L|_{C_x})$ (note that $d=1$).
By intersecting $X$ with a hyperplane containing $x$ and one of the lines,
we produce a $k$-rational conic passing through $x$, and approximating on the conic gives us
$2\geq \alpha_{x}(L)$.   We therefore conclude that $\alpha_{x}(L)=\min(2,\alpha_{x,C_x}(L|_{C_x}))$.

The curve $C_x$ is
singular at $x$, and since $x$ does not lie on a line, $C_x$ also
irreducible.  In particular, $C_x$ is an irreducible curve of
geometric genus zero, and since $x$ is defined over $k$, $C_x$ is
birational to $\P^1$ over $k$, via projection from $x$ in the tangent
plane.

Applying Theorem \ref{thm:curve} to $C_x$, we find that
$$\alpha_{x,C_x}(L|_{C_x}) =
\begin{cases}
\frac{3}{2} & \mbox{if $C_x$:}
\mbox{\begin{minipage}[t]{0.6\textwidth}
\begin{itemize}
\item[$\circ$] is cuspidal, or
\item[$\circ$] is nodal and the tangent lines have slopes in $k_v$ but not in $k$
\end{itemize}
\end{minipage}} \rule[-1.2cm]{0cm}{1.2cm}\\
3 & \mbox{if $C_x$ is nodal and the slopes of the tangent lines are in $k$}\\
\infty & \mbox{if $C_x$ is nodal and the slopes of the tangent lines are not in $k_v$} \\
\end{cases},
$$
and this implies the stated values of $\alpha_{x}(L)$ above. \epf

\vspace{.1in}

We now treat the case of a general nef divisor $D$.  In what follows,
we assume that the point $x$ does not lie on a $(-1)$-curve on $X$.  We
begin with a calculation of the Seshadri constant $\ep$.  To do
this, we will need some notation.

Let $\phi\colon X\to\P^2$ be the blowing down map, and let
$E_1,\ldots,E_6$ be the exceptional divisors of $\phi$.  We define the
following linear equivalence classes on $X$:
\begin{itemize}
\item $L=\phi^*\O(1)$
\item $L_i=L-E_i$, the strict transform of a line through $P_i=\phi(E_i)$
\item $L_{ij}=2L-(\sum E_n)+E_i+E_j$, the strict transform of a conic
through the four points $P_n$ with $n\neq i,j$
\item $B_i=3L-(\sum E_n)-E_i$, the strict transform of a cubic curve through
all six points $P_n$, with a node at $P_i$.
\end{itemize}

Let $\hyp$ be the class of a hyperplane in the anticanonical embedding $X\subset\PP^3$.
For any line $\ell$ on $X$, the hyperplanes
containing $\ell$ give (after removing $\ell$) a base-point-free pencil on $X$.
If $x\in X$ does not lie on a line then the unique curve in this pencil through $x$ is smooth and irreducible.
The classes $\{L_i, L_{ij}, B_i\}$ defined above
are the $27$ pencils coming from the lines.
Recall that for any point $x$ on $X$ we use $C_x$ for the
intersection of $X$ with its tangent plane at $x$ (so $C_x$ has class $\hyp$).
If $x$ does not lie on a line, then $C_x$ is a plane cubic curve with one double
point, at $x$.

\begin{theorem}\label{cubicseshadri}
Let $x$ be a point on $X$ that does not lie on a $(-1)$-curve, and let $D$ be
a nef divisor on $X$.  The Seshadri constant $\ep_x(D)$ is equal to
$\min\{D.L_i,D.L_{ij},D.B_i,(D.\hyp)/2\}$.
\end{theorem}

\noindent
{\it Proof:} \/ The nef cone $\Gamma$ of $X$ has $99$ generators, which
are listed in \S\ref{sec:appendix}, Table~\ref{nefgens}.
Let $S$ be the set of
$27$ divisor classes $\{L_i,L_{ij},B_i\}$ as $i$ and $j$ range over all
possible values, and for each element $C$ in $S$, we define the
subcone $\Gamma(C)$ by:
\[\Gamma(C)=\left\{D\in\Gamma\mid D.C=\min_{C'\in S}\{D.C'\}\,\mbox{and}\,
D.C\leq (D.\hyp)/2\right\}.\]
Further define the subcone $\Gamma(\hyp)$ to be:
\[\Gamma(\hyp)=\left\{D\in\Gamma\mid
(D.\hyp)/2\leq\min_{C'\in S}\{D.C'\}\right\}.\]
It is clear that $\Gamma$ is the union of these $28$ subcones.  To prove
Theorem~\ref{cubicseshadri}, it suffices to show that for every
subcone $\Gamma(C)$, with $C\in S$, the curve through $x$ in the pencil corresponding to $C$
is a Seshadri curve for $x$ with respect to $D$ for all $D\in \Gamma(C)$
(repectively, in the case of the subcone $\Gamma(\hyp)$, that  $C_x$ is a Seshadri curve for $x$ with respect
to $D$ for all $D\in \Gamma(\hyp)$).
By Proposition~\ref{prop:alphacurve}(b) it further
suffices to prove this for $D$ a generator of the cone $\Gamma(C)$ (respectively $\Gamma(\hyp)$).

The fundamental group of the space of all smooth cubic surfaces acts via monodromy on the
N\'eron-Severi lattice of $X$.  This monodromy action preserves the hyperplane class $h$ and acts transitively
on the classes of the lines.  Thus, up to monodromy action, there are only two of these subcones:
$\Gamma(L_1)$ and $\Gamma(\hyp)$.
Generators for each of these
subcones can also be found in \S\ref{sec:appendix}.  Let $F=F_{x,L_1}$ be the unique curve in the pencil $L_1$
passing through $x$. For
each generator $D$ of $\Gamma(L_1)$, it is straightforward to verify
that $F$ is a Seshadri curve for $x$ with respect to $D$. These verifications
also appear in \S\ref{sec:appendix}.  Each
generator $G$ of $\Gamma(\hyp)$ is also a generator of one of the other
twenty-seven subcones $\Gamma(C)$, and for each such $G$, we have
$G.C=(G.\hyp)/2=(G.C_x)/\mult_{x}C_x$.  Thus, since $C$ is a Seshadri curve for $x$ with respect to $G$, it
follows that $C_x$ is also a Seshadri curve for $x$ with respect to $G$, and so $C_x$ is a
Seshadri curve for every element of $\Gamma(\hyp)$.  This concludes the
proof.  \qed

\vspace{.1in}

The next step is to calculate $\alpha_x$ for a point on a cubic surface.

\begin{theorem}\label{cubicalpha}
Let $x$ be a point on $X$ that does not lie on a $(-1)$-curve, and let
$D$ be a nef divisor on $X$.  If the tangent curve $C_x$ is a cuspidal
cubic, or a nodal cubic whose tangent lines at $x$ are defined over
$k_v$ but not defined over $k$, then $\alpha_x(D)=\ep_x(D)$.  Otherwise,
$\alpha_x(D)=\min\{D.L_i,D.L_{ij},D.B_i\}$.
\end{theorem}

\noindent

{\it Proof:}
Suppose that $D$ is in one of the cones $\Gamma(C)$ for $C\in S$, and let $F_{x,C}$ be the element of
the pencil corresponding to $C$ which passes through $x$. Since $F_{x,C}$ is a smooth $k$-rational curve,
we have
$$D.C \stackrel{\scriptsize\ref{thm:curve}}{=} \alpha_{x}(D|_{F_{x,C}})
\stackrel{\scriptsize\ref{prop:alpha-and-ep}(c)}{\geq} \alpha_{x}(D)
\stackrel{\scriptsize\ref{cor:alpha-vs-epsilon}}{\geq} \ep_{x}(D)
\stackrel{\scriptsize\ref{thm:cubic}}{=} D.C,$$
where, reading from left to right, the equalities and inequalites are given by
Theorem \ref{thm:curve}, Proposition \ref{prop:alpha-and-ep}(c), Corollary \ref{cor:alpha-vs-epsilon},
and Theorem \ref{thm:cubic} respectively.   Thus $\alpha_{x}(D) = D.C$ and $F_{x,C}$ is a curve of best
approximation with respect to $D$.

Now suppose that $D\in \Gamma(\hyp)$.  If $C_{x}$ is cuspidal, or nodal with tangent lines having slopes
in $k_v$ but not $k$, then Theorem \ref{thm:curve} gives $\alpha_{x}(D|_{C_x})=D.C_x/2=D.C_x/\mult_{x}C_x$.
By Theorem \ref{thm:cubic}, $\ep_{x}(D)=D.C_{x}/2$, and so as above we conclude that
$\alpha_{x}(D)=D.C_x=\ep_{x}(D)$, and that $C_x$ is a curve of best approximation for $D$.

We now assume that $C_{x}$ is nodal and the slopes of the tangent lines are in $k$ or not in $k_v$.
The codimension one faces of $\Gamma(\hyp)$ (i.e., the facets) occur where one of the inequalities 
defining $\Gamma(\hyp)$ becomes an equality, so that each facet is the intersection of $\Gamma(\hyp)$ and
$\Gamma(C)$ for some $C\in S$.  For each $C\in S$ set $\hat{\Gamma}(C)$ to be the cone generated by 
$\Gamma(C)$ and $-K$.  Since $-K$ is in the interior of $\Gamma(\hyp)$ it follows that $\Gamma$ is the union
of the $\hat{\Gamma}(C)$, $C\in S$.

For any $C\in S$, let $F_{x,C}$ be the member of the pencil corresponding to $C$ passing through $x$, as
in the first part of the argument.
In the proof of Theorem \ref{thm:cubic} we have seen that $F_{x,C}$ is a curve of best approximation for $-K$,
and in the first part of the argument above that $F_{x,C}$ is a curve of best approximation for all $D\in \Gamma(C)$.
By Proposition \ref{prop:alphacurve}(a) we conclude that $F_{x,C}$ is a curve of best approximation for all
$D\in \hat{\Gamma}(C)$.  The result follows.  \qed

\vspace{.1in}

Note that as part of the proof we have shown that Conjecture~\ref{conj:ratcurve} holds for every point $x\in X$ not on a $(-1)$-curve.

\section{Appendix: Generators of nef cones and subcones for the cubic surface}
\label{sec:appendix}

A version of this appendix, with additional tables and larger font, may be found at
\cite{Porta}.
We use the notation from \S\ref{sec:cubic}.
In each of the tables in this appendix the first column is a numerical
identifier of the vector in that row.  The subsequent columns
represent the coefficients of the vector with respect to the basis
$\{L,E_1,\ldots,E_6\}$ of the N\'eron-Severi group of $X$.  Thus,
vector number 1 in Table~\ref{nefgens} is the divisor class
$2L-E_1-E_2-E_3$.
Each of the cones has $99$ generators.
There is no correspondence or relation between rows in different tables with the same numerical identifier.

Table~\ref{nefgens}, of generators of the nef cone, is reproducing
information that has been well known for some time, of course.  
It was
calculated for these tables by finding generators for the cone obtained
as the intersection of the half-spaces corresponding to non-negative
intersection with each of the 27 lines on the cubic surface.
The other tables were generated in a similar way.  For instance,
Table~\ref{ligens}, of generators of the cone $\Gamma(L_1)$, was generated by
using the half-spaces defining $\Gamma$, in addition to the half-spaces
corresponding to the intersection inequalities described above.

{
\begin{centering}
\tablefont
\begin{tabular}{|c|c|c|}
\multicolumn{3}{c}{\bf \normalsize Table \tabnum: Generators of the nef cone $\Gamma$ of a smooth cubic surface \label{nefgens} } \\
\multicolumn{3}{c}{} \\
\hline
\begin{tabular}{lrrrrrrr}
\bf \# & $L$ & $E_1$ & $E_2$ & $E_3$ & $E_4$ & $E_5$ & $E_6$ \\
\hline
\bf  1 & 2 & -1 & -1 & -1 &  0 &  0 &  0 \\
\bf   2 & 2 & -1 & -1 &  0 & -1 &  0 &  0 \\
\bf   3 & 2 & -1 & -1 &  0 &  0 & -1 &  0 \\
\bf   4  & 2  & -1  & -1  &  0  &  0  &  0  & -1  \\
\bf   5  & 2  & -1  &  0  & -1  & -1  &  0  &  0  \\
\bf   6  & 2  & -1  &  0  & -1  &  0  & -1  &  0  \\
\bf   7  & 2  & -1  &  0  & -1  &  0  &  0  & -1  \\
\bf   8  & 2  & -1  &  0  &  0  & -1  & -1  &  0  \\
\bf   9  & 2  & -1  &  0  &  0  & -1  &  0  & -1  \\
\bf  10  & 2  & -1  &  0  &  0  &  0  & -1  & -1  \\
\bf  11  & 1  &  0  &  0  &  0  &  0  &  0  &  0  \\
\bf  12  & 3  & -2  & -1  & -1  & -1  & -1  &  0  \\
\bf  13  & 3  & -2  & -1  & -1  & -1  &  0  & -1  \\
\bf  14  & 3  & -2  & -1  & -1  &  0  & -1  & -1  \\
\bf  15  & 3  & -2  & -1  &  0  & -1  & -1  & -1  \\
\bf  16  & 3  & -2  &  0  & -1  & -1  & -1  & -1  \\
\bf  17  & 1  & -1  &  0  &  0  &  0  &  0  &  0  \\
\bf  18  & 1  &  0  & -1  &  0  &  0  &  0  &  0  \\
\bf  19  & 1  &  0  &  0  & -1  &  0  &  0  &  0  \\
\bf  20  & 1  &  0  &  0  &  0  & -1  &  0  &  0  \\
\bf  21  & 1  &  0  &  0  &  0  &  0  & -1  &  0  \\
\bf  22  & 1  &  0  &  0  &  0  &  0  &  0  & -1  \\
\bf  23  & 2  &  0  & -1  & -1  & -1  &  0  &  0  \\
\bf  24  & 2  &  0  & -1  & -1  &  0  & -1  &  0  \\
\bf  25  & 2  &  0  & -1  & -1  &  0  &  0  & -1  \\
\bf  26  & 2  &  0  & -1  &  0  & -1  & -1  &  0  \\
\bf  27  & 2  &  0  & -1  &  0  & -1  &  0  & -1  \\
\bf  28  & 2  &  0  & -1  &  0  &  0  & -1  & -1  \\
\bf  29  & 2  &  0  &  0  & -1  & -1  & -1  &  0  \\
\bf  30  & 2  &  0  &  0  & -1  & -1  &  0  & -1  \\
\bf  31  & 2  &  0  &  0  & -1  &  0  & -1  & -1  \\
\bf  32  & 2  &  0  &  0  &  0  & -1  & -1  & -1  \\
\bf  33  & 2  & -1  & -1  & -1  & -1  &  0  &  0  \\
\end{tabular}
&
\begin{tabular}{lrrrrrrr}
\bf \# & $L$ & $E_1$ & $E_2$ & $E_3$ & $E_4$ & $E_5$ & $E_6$ \\
\hline
\bf  34  & 2  & -1  & -1  & -1  &  0  & -1  &  0  \\
\bf  35  & 2  & -1  & -1  & -1  &  0  &  0  & -1  \\
\bf  36  & 2  & -1  & -1  &  0  & -1  & -1  &  0  \\
\bf  37  & 2  & -1  & -1  &  0  & -1  &  0  & -1  \\
\bf  38  & 2  & -1  & -1  &  0  &  0  & -1  & -1  \\
\bf  39  & 2  & -1  &  0  & -1  & -1  & -1  &  0  \\
\bf  40  & 2  & -1  &  0  & -1  & -1  &  0  & -1  \\
\bf  41  & 2  & -1  &  0  & -1  &  0  & -1  & -1  \\
\bf  42  & 2  & -1  &  0  &  0  & -1  & -1  & -1  \\
\bf  43  & 2  &  0  & -1  & -1  & -1  & -1  &  0  \\
\bf  44  & 2  &  0  & -1  & -1  & -1  &  0  & -1  \\
\bf  45  & 2  &  0  & -1  & -1  &  0  & -1  & -1  \\
\bf  46  & 2  &  0  & -1  &  0  & -1  & -1  & -1  \\
\bf  47  & 2  &  0  &  0  & -1  & -1  & -1  & -1  \\
\bf  48  & 3  & -1  & -2  & -1  & -1  & -1  &  0  \\
\bf  49  & 3  & -1  & -2  & -1  & -1  &  0  & -1  \\
\bf  50  & 3  & -1  & -2  & -1  &  0  & -1  & -1  \\
\bf  51  & 3  & -1  & -2  &  0  & -1  & -1  & -1  \\
\bf  52  & 3  & -1  & -1  & -2  & -1  & -1  &  0  \\
\bf  53  & 3  & -1  & -1  & -2  & -1  &  0  & -1  \\
\bf  54  & 3  & -1  & -1  & -2  &  0  & -1  & -1  \\
\bf  55  & 3  & -1  & -1  & -1  & -2  & -1  &  0  \\
\bf  56  & 3  & -1  & -1  & -1  & -2  &  0  & -1  \\
\bf  57  & 3  & -1  & -1  & -1  & -1  & -2  &  0  \\
\bf  58  & 3  & -1  & -1  & -1  & -1  &  0  & -2  \\
\bf  59  & 3  & -1  & -1  & -1  &  0  & -2  & -1  \\
\bf  60  & 3  & -1  & -1  & -1  &  0  & -1  & -2  \\
\bf  61  & 3  & -1  & -1  &  0  & -2  & -1  & -1  \\
\bf  62  & 3  & -1  & -1  &  0  & -1  & -2  & -1  \\
\bf  63  & 3  & -1  & -1  &  0  & -1  & -1  & -2  \\
\bf  64  & 3  & -1  &  0  & -2  & -1  & -1  & -1  \\
\bf  65  & 3  & -1  &  0  & -1  & -2  & -1  & -1  \\
\bf  66  & 3  & -1  &  0  & -1  & -1  & -2  & -1  \\
\end{tabular}
&
\begin{tabular}{lrrrrrrr}
\bf \# & $L$ & $E_1$ & $E_2$ & $E_3$ & $E_4$ & $E_5$ & $E_6$ \\
\hline
\bf  67  & 3  & -1  &  0  & -1  & -1  & -1  & -2  \\
\bf  68  & 3  &  0  & -2  & -1  & -1  & -1  & -1  \\
\bf  69  & 3  &  0  & -1  & -2  & -1  & -1  & -1  \\
\bf  70  & 3  &  0  & -1  & -1  & -2  & -1  & -1  \\
\bf  71  & 3  &  0  & -1  & -1  & -1  & -2  & -1  \\
\bf  72  & 3  &  0  & -1  & -1  & -1  & -1  & -2  \\
\bf  73  & 3  & -2  & -1  & -1  & -1  & -1  & -1  \\
\bf  74  & 3  & -1  & -2  & -1  & -1  & -1  & -1  \\
\bf  75  & 3  & -1  & -1  & -2  & -1  & -1  & -1  \\
\bf  76  & 3  & -1  & -1  & -1  & -2  & -1  & -1  \\
\bf  77  & 3  & -1  & -1  & -1  & -1  & -2  & -1  \\
\bf  78  & 3  & -1  & -1  & -1  & -1  & -1  & -2  \\
\bf  79  & 4  & -2  & -2  & -2  & -1  & -1  & -1  \\
\bf  80  & 4  & -2  & -2  & -1  & -2  & -1  & -1  \\
\bf  81  & 4  & -2  & -2  & -1  & -1  & -2  & -1  \\
\bf  82  & 4  & -2  & -2  & -1  & -1  & -1  & -2  \\
\bf  83  & 4  & -2  & -1  & -2  & -2  & -1  & -1  \\
\bf  84  & 4  & -2  & -1  & -2  & -1  & -2  & -1  \\
\bf  85  & 4  & -2  & -1  & -2  & -1  & -1  & -2  \\
\bf  86  & 4  & -2  & -1  & -1  & -2  & -2  & -1  \\
\bf  87  & 4  & -2  & -1  & -1  & -2  & -1  & -2  \\
\bf  88  & 4  & -2  & -1  & -1  & -1  & -2  & -2  \\
\bf  89  & 4  & -1  & -2  & -2  & -2  & -1  & -1  \\
\bf  90  & 4  & -1  & -2  & -2  & -1  & -2  & -1  \\
\bf  91  & 4  & -1  & -2  & -2  & -1  & -1  & -2  \\
\bf  92  & 4  & -1  & -2  & -1  & -2  & -2  & -1  \\
\bf  93  & 4  & -1  & -2  & -1  & -2  & -1  & -2  \\
\bf  94  & 4  & -1  & -2  & -1  & -1  & -2  & -2  \\
\bf  95  & 4  & -1  & -1  & -2  & -2  & -2  & -1  \\
\bf  96  & 4  & -1  & -1  & -2  & -2  & -1  & -2  \\
\bf  97  & 4  & -1  & -1  & -2  & -1  & -2  & -2  \\
\bf  98  & 4  & -1  & -1  & -1  & -2  & -2  & -2  \\
\bf  99  & 5  & -2  & -2  & -2  & -2  & -2  & -2
\end{tabular} \\
\hline
\end{tabular}\\
\end{centering}
}

\bigskip
\bigskip
\bigskip
In Table~\ref{ligens} which follows, we use $D_n$ to refer to the divisor class
represented by row $n$ of Table~\ref{ligens}.
For any point $x\in X$ not on a $(-1)$-curve, the unique curve $F=F_{x,L_1}$ in the pencil $L_1$ passing through $x$ is smooth and irreducible.
In each line of the table ``Reason'' is a --- very brief! --- justification of why $F$ is a
Seshadri curve for $x$ with respect to $D_n$.

For instance, in row 1 of Table~\ref{ligens}, the ``Reason'' is $L_1.D_1=1$, and thus
$F.D_1=L_1.D_1=1$.  
We claim that for the divisors $D_i$, $\epsilon_x$ is always at least one if it is nonzero.  To see this, notice 
that the generators of the nef cone (see Table~\ref{nefgens}) are all either morphisms to $\P^1$ corresponding to 
pencils of conics on the cubic surface, or else morphisms to $\P^2$ that are the blowing down of six pairwise disjoint
$(-1)$-curves.  In both cases, the Seshadri constant is easily seen to be either zero or at least one.  It is 
straightforward to check that all the generators listed in Table~\ref{ligens} are non-negative integer linear 
combinations of the generators of the nef cones, and therefore (by Proposition \ref{prop:alpha-and-ep}(b))
enjoy the same property: for any point $x$, 
the Seshadri constant $\epsilon_{x}(D_i)$ is either zero or else is at least one.

By assumption, $x$ does not lie on any $(-1)$-curve, which are the only curves contracted by any $D_i$ (except for $D_{18}=L_1$, for which $\epsilon=0$ for all points).  Therefore, since $F$ has degree $1$ with respect to $D_1$, $F$ 
is a Seshadri curve for $x$ with respect to $D_1$.

As a second example, 
in row 29 of Table~\ref{ligens}, the comment ``$L+L_{56}$'' means that the divisor $D_{29}$ represented by that row is the sum of $L$ and $L_{56}$.  Any curve that has nonzero intersection with $L$ must have $L.C/\mbox{mult}_x(C)\geq 1$, for any $x$ not lying on a $(-1)$-curve, since $L$ is an isomorphism away from $(-1)$-curves.  Similarly, any curve not contracted by $L_{56}$ must also satisfy $L_{56}.C/\mbox{mult}_x(C)\geq 1$, so any curve not contracted by $L_{56}$ or $L$ must satisfy $(L+L_{56}).C/\mbox{mult}_x(C)\geq 2$.  If $C$ is contracted by $L_{56}$, then it is either a $(-1)$-curve, or else it is an element of the divisor class $L_{56}$ itself, in which case it satisfies $(L+L_{56}).C/\mbox{mult}_x(C)=2$ by direct calculation.  In all cases, since $x$ does not lie on a $(-1)$-curve, we see that $\epsilon_{x}(L+L_{56})\geq2$, and since $L_1.L=L_1.L_{56}=1$, we conclude that $\epsilon_x(L+L_{56})=2$, and so the curve in the class $L_1$ through $x$ is a Seshadri curve for $x$ with respect to $D_{29}=L+L_{56}$.  Similar arguments explain the 
other reasons of the form ``$A+B$'' or ``$A+B+C$''.

In light of these arguments, for Table~\ref{ligens}, it is useful to know that $L_1$ has intersection number one with 
the divisors $L$, $B_1$, $L_i$ for $i\neq 1$, and $L_{ij}$ for $i,j\neq 1$.

\begin{centering}
\tablefont
\begin{tabular}{|c|c|}
\multicolumn{2}{c}{\bf \normalsize Table \tabnum: Generators of the cone $\Gamma(L_1)$ \label{ligens}} \\
\multicolumn{2}{c}{} \\
\hline
\begin{tabular}{rrrrrrrrr}
\bf \# & $L$ & $E_1$ & $E_2$ & $E_3$ & $E_4$ & $E_5$ & $E_6$ & Reason \\
\hline
\bf  1  & 4  & -3  & -1  & -1  & -1  & -1  & -1  & $L_1.D_1=1$ \\
\bf  2  & 2  & -1  & -1  &  0  &  0  &  0  &  0  & $L_1.D_2=1$ \\
\bf  3  & 2  & -1  &  0  & -1  &  0  &  0  &  0  & $L_1.D_3=1$ \\
\bf  4  & 2  & -1  &  0  &  0  & -1  &  0  &  0  & $L_1.D_4=1$ \\
\bf  5  & 2  & -1  &  0  &  0  &  0  & -1  &  0  & $L_1.D_5=1$ \\
\bf  6  & 2  & -1  &  0  &  0  &  0  &  0  & -1  & $L_1.D_6=1$ \\
\bf  7  & 1  &  0  &  0  &  0  &  0  &  0  &  0  & $L_1.D_7=1$ \\
\bf  8  & 3  & -2  & -1  & -1  & -1  &  0  &  0  & $L_1.D_8=1$ \\
\bf  9  & 3  & -2  & -1  & -1  &  0  & -1  &  0  & $L_1.D_9=1$ \\
\bf 10  & 3  & -2  & -1  & -1  &  0  &  0  & -1  & $L_1.D_{10}=1$ \\
\bf 11  & 3  & -2  & -1  &  0  & -1  & -1  &  0  & $L_1.D_{11}=1$ \\
\bf 12  & 3  & -2  & -1  &  0  & -1  &  0  & -1  & $L_1.D_{12}=1$ \\
\bf 13  & 3  & -2  & -1  &  0  &  0  & -1  & -1  & $L_1.D_{13}=1$ \\
\bf 14  & 3  & -2  &  0  & -1  & -1  & -1  &  0  & $L_1.D_{14}=1$ \\
\bf 15  & 3  & -2  &  0  & -1  & -1  &  0  & -1  & $L_1.D_{15}=1$ \\
\bf 16  & 3  & -2  &  0  & -1  &  0  & -1  & -1  & $L_1.D_{16}=1$ \\
\bf 17  & 3  & -2  &  0  &  0  & -1  & -1  & -1  & $L_1.D_{17}=1$ \\
\bf 18  & 1  & -1  &  0  &  0  &  0  &  0  &  0  & $L_1.D_{18}=0$ \\
\bf 19  & 2  & -1  & -1  & -1  &  0  &  0  &  0  & $L_1.D_{19}=1$ \\
\bf 20  & 2  & -1  & -1  &  0  & -1  &  0  &  0  & $L_1.D_{20}=1$ \\
\bf 21  & 2  & -1  & -1  &  0  &  0  & -1  &  0  & $L_1.D_{21}=1$ \\
\bf 22  & 2  & -1  & -1  &  0  &  0  &  0  & -1  & $L_1.D_{22}=1$ \\
\bf 23  & 2  & -1  &  0  & -1  & -1  &  0  &  0  & $L_1.D_{23}=1$ \\
\bf 24  & 2  & -1  &  0  & -1  &  0  & -1  &  0  & $L_1.D_{24}=1$ \\
\bf 25  & 2  & -1  &  0  & -1  &  0  &  0  & -1  & $L_1.D_{25}=1$ \\
\bf 26  & 2  & -1  &  0  &  0  & -1  & -1  &  0  & $L_1.D_{26}=1$ \\
\bf 27  & 2  & -1  &  0  &  0  & -1  &  0  & -1  & $L_1.D_{27}=1$ \\
\bf 28  & 2  & -1  &  0  &  0  &  0  & -1  & -1  & $L_1.D_{28}=1$ \\
\bf 29  & 3  & -1  & -1  & -1  & -1  &  0  &  0  & $L+L_{56}$ \\
\bf 30  & 3  & -1  & -1  & -1  &  0  & -1  &  0  & $L+L_{46}$ \\
\bf 31  & 3  & -1  & -1  & -1  &  0  &  0  & -1  & $L+L_{45}$ \\
\bf 32  & 3  & -1  & -1  &  0  & -1  & -1  &  0  & $L+L_{36}$ \\
\bf 33  & 3  & -1  & -1  &  0  & -1  &  0  & -1  & $L+L_{35}$ \\
\bf 34  & 3  & -1  & -1  &  0  &  0  & -1  & -1  & $L+L_{34}$ \\
\bf 35  & 3  & -1  &  0  & -1  & -1  & -1  &  0  & $L+L_{26}$ \\
\bf 36  & 3  & -1  &  0  & -1  & -1  &  0  & -1  & $L+L_{25}$ \\
\bf 37  & 3  & -1  &  0  & -1  &  0  & -1  & -1  & $L+L_{24}$ \\
\bf 38  & 3  & -1  &  0  &  0  & -1  & -1  & -1  & $L+L_{23}$ \\
\bf 39  & 3  & -1  & -1  & -1  & -1  & -1  &  0  & $L_2+L_{26}$ \\
\bf 40  & 3  & -1  & -1  & -1  & -1  &  0  & -1  & $L_2+L_{25}$ \\
\bf 41  & 3  & -1  & -1  & -1  &  0  & -1  & -1  & $L_2+L_{24}$ \\
\bf 42  & 3  & -1  & -1  &  0  & -1  & -1  & -1  & $L_2+L_{23}$ \\
\bf 43  & 3  & -1  &  0  & -1  & -1  & -1  & -1  & $L_3+L_{23}$ \\
\bf 44  & 4  & -1  & -1  & -1  & -1  & -1  & -1  & $L_{23}+L_2+L_3$ \\
\bf 45  & 3  & -2  & -1  & -1  & -1  & -1  &  0  & $L_1.D_{45}=1$ \\
\bf 46  & 3  & -2  & -1  & -1  & -1  &  0  & -1  & $L_1.D_{46}=1$ \\
\bf 47  & 3  & -2  & -1  & -1  &  0  & -1  & -1  & $L_1.D_{47}=1$ \\
\bf 48  & 3  & -2  & -1  &  0  & -1  & -1  & -1  & $L_1.D_{48}=1$ \\
\bf 49  & 3  & -2  &  0  & -1  & -1  & -1  & -1  & $L_1.D_{49}=1$ \\
\bf 50  & 4  & -2  & -2  & -1  & -1  & -1  &  0  & $D_{45}+L_2$ \\
\end{tabular}
&
\begin{tabular}{rrrrrrrrr}
\bf \# & $L$ & $E_1$ & $E_2$ & $E_3$ & $E_4$ & $E_5$ & $E_6$ & Reason \\
\hline
\bf 51  & 4  & -2  & -2  & -1  & -1  &  0  & -1  & $D_{46}+L_2$ \\
\bf 52  & 4  & -2  & -2  & -1  &  0  & -1  & -1  & $D_{47}+L_2$ \\
\bf 53  & 4  & -2  & -2  &  0  & -1  & -1  & -1  & $D_{48}+L_2$ \\
\bf 54  & 4  & -2  & -1  & -2  & -1  & -1  &  0  & $D_{45}+L_3$ \\
\bf 55  & 4  & -2  & -1  & -2  & -1  &  0  & -1  & $D_{46}+L_3$ \\
\bf 56  & 4  & -2  & -1  & -2  &  0  & -1  & -1  & $D_{47}+L_3$ \\
\bf 57  & 4  & -2  & -1  & -1  & -2  & -1  &  0  & $D_{45}+L_4$ \\
\bf 58  & 4  & -2  & -1  & -1  & -2  &  0  & -1  & $D_{46}+L_4$ \\
\bf 59  & 4  & -2  & -1  & -1  & -1  & -2  &  0  & $D_{45}+L_5$ \\
\bf 60  & 4  & -2  & -1  & -1  & -1  &  0  & -2  & $D_{46}+L_6$ \\
\bf 61  & 4  & -2  & -1  & -1  &  0  & -2  & -1  & $D_{47}+L_5$ \\
\bf 62  & 4  & -2  & -1  & -1  &  0  & -1  & -2  & $D_{47}+L_6$ \\
\bf 63  & 4  & -2  & -1  &  0  & -2  & -1  & -1  & $D_{48}+L_4$ \\
\bf 64  & 4  & -2  & -1  &  0  & -1  & -2  & -1  & $D_{48}+L_5$ \\
\bf 65  & 4  & -2  & -1  &  0  & -1  & -1  & -2  & $D_{48}+L_6$ \\
\bf 66  & 4  & -2  &  0  & -2  & -1  & -1  & -1  & $D_{49}+L_3$ \\
\bf 67  & 4  & -2  &  0  & -1  & -2  & -1  & -1  & $D_{49}+L_4$ \\
\bf 68  & 4  & -2  &  0  & -1  & -1  & -2  & -1  & $D_{49}+L_5$ \\
\bf 69  & 4  & -2  &  0  & -1  & -1  & -1  & -2  & $D_{49}+L_6$ \\
\bf 70  & 4  & -2  & -2  & -1  & -1  & -1  & -1  & $B_1+L_2$ \\
\bf 71  & 4  & -2  & -1  & -2  & -1  & -1  & -1  & $B_1+L_3$ \\
\bf 72  & 4  & -2  & -1  & -1  & -2  & -1  & -1  & $B_1+L_4$ \\
\bf 73  & 4  & -2  & -1  & -1  & -1  & -2  & -1  & $B_1+L_5$ \\
\bf 74  & 4  & -2  & -1  & -1  & -1  & -1  & -2  & $B_1+L_6$ \\
\bf 75  & 5  & -2  & -2  & -2  & -1  & -1  & -1  & $L_3+L_{34}+L_{56}$ \\
\bf 76  & 5  & -2  & -2  & -1  & -2  & -1  & -1  & $L_2+L_{23}+L_{56}$ \\
\bf 77  & 5  & -2  & -2  & -1  & -1  & -2  & -1  & $L_2+L_{23}+L_{46}$ \\
\bf 78  & 5  & -2  & -2  & -1  & -1  & -1  & -2  & $L_2+L_{23}+L_{45}$ \\
\bf 79  & 5  & -2  & -1  & -2  & -2  & -1  & -1  & $L_3+L_{23}+L_{56}$ \\
\bf 80  & 5  & -2  & -1  & -2  & -1  & -2  & -1  & $L_3+L_{23}+L_{46}$ \\
\bf 81  & 5  & -2  & -1  & -2  & -1  & -1  & -2  & $L_3+L_{23}+L_{45}$ \\
\bf 82  & 5  & -2  & -1  & -1  & -2  & -2  & -1  & $L_4+L_{23}+L_{46}$ \\
\bf 83  & 5  & -2  & -1  & -1  & -2  & -1  & -2  & $L_4+L_{23}+L_{45}$ \\
\bf 84  & 5  & -2  & -1  & -1  & -1  & -2  & -2  & $L_5+L_{23}+L_{45}$ \\
\bf 85  & 5  & -3  & -2  & -2  & -1  & -1  & -1  & $D_{45}+L_{45}$ \\
\bf 86  & 5  & -3  & -2  & -1  & -2  & -1  & -1  & $D_{45}+L_{35}$ \\
\bf 87  & 5  & -3  & -2  & -1  & -1  & -2  & -1  & $D_{45}+L_{34}$ \\
\bf 88  & 5  & -3  & -2  & -1  & -1  & -1  & -2  & $D_{46}+L_{34}$ \\
\bf 89  & 5  & -3  & -1  & -2  & -2  & -1  & -1  & $D_{45}+L_{25}$ \\
\bf 90  & 5  & -3  & -1  & -2  & -1  & -2  & -1  & $D_{45}+L_{24}$ \\
\bf 91  & 5  & -3  & -1  & -2  & -1  & -1  & -2  & $D_{46}+L_{24}$ \\
\bf 92  & 5  & -3  & -1  & -1  & -2  & -2  & -1  & $D_{45}+L_{23}$ \\
\bf 93  & 5  & -3  & -1  & -1  & -2  & -1  & -2  & $D_{46}+L_{23}$ \\
\bf 94  & 5  & -3  & -1  & -1  & -1  & -2  & -2  & $D_{47}+L_{23}$ \\
\bf 95  & 6  & -3  & -2  & -2  & -2  & -2  & -1  & $L_{23}+L_{46}+L_{56}$ \\
\bf 96  & 6  & -3  & -2  & -2  & -2  & -1  & -2  & $L_{23}+L_{45}+L_{56}$ \\
\bf 97  & 6  & -3  & -2  & -2  & -1  & -2  & -2  & $L_{23}+L_{45}+L_{46}$ \\
\bf 98  & 6  & -3  & -2  & -1  & -2  & -2  & -2  & $L_{23}+L_{34}+L_{56}$ \\
\bf 99  & 6  & -3  & -1  & -2  & -2  & -2  & -2  & $L_{23}+L_{24}+L_{56}$ \\
\\
\end{tabular} \\
\hline
\end{tabular} \\
\end{centering}

\vfill

In Table~\ref{cpgens}, the rightmost column of row $n$ contains a
divisor class $C\in S$ such that $G_n$
(the divisor corresponding to the $n$th row of
Table \ref{cpgens}) is also a generator of the subcone $\Gamma(C)$.
From the definition of the cones $\Gamma(C)$ and $\Gamma(\hyp)$, this implies that
$G_n.C=(G_n.\hyp)/2$. As explained in the proof of Theorem \ref{cubicalpha},
this provides a verification that $C_x$ is a Seshadri curve for $x$ with respect to $G_n$.

\vfill
\phantom{.}

\begin{centering}
\tablefont
\begin{tabular}{|c|c|}
\multicolumn{2}{c}{\bf\normalsize Table \tabnum: Generators of the cone $\Gamma(\hyp)$ \label{cpgens}} \\
\multicolumn{2}{c}{} \\
\hline
\begin{tabular}{rrrrrrrrr}
\bf \# & $L$ & $E_1$ & $E_2$ & $E_3$ & $E_4$ & $E_5$ & $E_6$ & Divisor Class \\
\hline
\bf  1  & 8  & -3  & -3  & -3  & -3  & -3  & -3  & $B_1$ \\
\bf  2  & 4  & -1  & -1  & -1  & -1  & -1  & -1  & $L_1$ \\
\bf  3  & 4  & -2  & -2  & -1  & -1  & -1  & -1  & $L_1$ \\
\bf  4  & 4  & -2  & -1  & -2  & -1  & -1  & -1  & $L_1$ \\
\bf  5  & 4  & -2  & -1  & -1  & -2  & -1  & -1  & $L_1$ \\
\bf  6  & 4  & -2  & -1  & -1  & -1  & -2  & -1  & $L_1$ \\
\bf  7  & 4  & -2  & -1  & -1  & -1  & -1  & -2  & $L_1$ \\
\bf  8  & 4  & -1  & -2  & -2  & -1  & -1  & -1  & $L_2$ \\
\bf  9  & 4  & -1  & -2  & -1  & -2  & -1  & -1  & $L_2$ \\
\bf 10  & 4  & -1  & -2  & -1  & -1  & -2  & -1  & $L_2$ \\
\bf 11  & 4  & -1  & -2  & -1  & -1  & -1  & -2  & $L_2$ \\
\bf 12  & 4  & -1  & -1  & -2  & -2  & -1  & -1  & $L_3$ \\
\bf 13  & 4  & -1  & -1  & -2  & -1  & -2  & -1  & $L_3$ \\
\bf 14  & 4  & -1  & -1  & -2  & -1  & -1  & -2  & $L_3$ \\
\bf 15  & 4  & -1  & -1  & -1  & -2  & -2  & -1  & $L_4$ \\
\bf 16  & 4  & -1  & -1  & -1  & -2  & -1  & -2  & $L_4$ \\
\bf 17  & 4  & -1  & -1  & -1  & -1  & -2  & -2  & $L_5$ \\
\bf 18  & 5  & -2  & -2  & -2  & -1  & -1  & -1  & $L_1$ \\
\bf 19  & 5  & -2  & -2  & -1  & -2  & -1  & -1  & $L_1$ \\
\bf 20  & 5  & -2  & -2  & -1  & -1  & -2  & -1  & $L_1$ \\
\bf 21  & 5  & -2  & -2  & -1  & -1  & -1  & -2  & $L_1$ \\
\bf 22  & 5  & -2  & -1  & -2  & -2  & -1  & -1  & $L_1$ \\
\bf 23  & 5  & -2  & -1  & -2  & -1  & -2  & -1  & $L_1$ \\
\bf 24  & 5  & -2  & -1  & -2  & -1  & -1  & -2  & $L_1$ \\
\bf 25  & 5  & -2  & -1  & -1  & -2  & -2  & -1  & $L_1$ \\
\bf 26  & 5  & -2  & -1  & -1  & -2  & -1  & -2  & $L_1$ \\
\bf 27  & 5  & -2  & -1  & -1  & -1  & -2  & -2  & $L_1$ \\
\bf 28  & 5  & -1  & -2  & -2  & -2  & -1  & -1  & $L_2$ \\
\bf 29  & 5  & -1  & -2  & -2  & -1  & -2  & -1  & $L_2$ \\
\bf 30  & 5  & -1  & -2  & -2  & -1  & -1  & -2  & $L_2$ \\
\bf 31  & 5  & -1  & -2  & -1  & -2  & -2  & -1  & $L_2$ \\
\bf 32  & 5  & -1  & -2  & -1  & -2  & -1  & -2  & $L_2$ \\
\bf 33  & 5  & -1  & -2  & -1  & -1  & -2  & -2  & $L_2$ \\
\bf 34  & 5  & -1  & -1  & -2  & -2  & -2  & -1  & $L_3$ \\
\bf 35  & 5  & -1  & -1  & -2  & -2  & -1  & -2  & $L_3$ \\
\bf 36  & 5  & -1  & -1  & -2  & -1  & -2  & -2  & $L_3$ \\
\bf 37  & 5  & -1  & -1  & -1  & -2  & -2  & -2  & $L_4$ \\
\bf 38  & 7  & -3  & -3  & -3  & -2  & -2  & -2  & $B_1$ \\
\bf 39  & 7  & -3  & -3  & -2  & -3  & -2  & -2  & $B_1$ \\
\bf 40  & 7  & -3  & -3  & -2  & -2  & -3  & -2  & $B_1$ \\
\bf 41  & 7  & -3  & -3  & -2  & -2  & -2  & -3  & $B_1$ \\
\bf 42  & 7  & -3  & -2  & -3  & -3  & -2  & -2  & $B_1$ \\
\bf 43  & 7  & -3  & -2  & -3  & -2  & -3  & -2  & $B_1$ \\
\bf 44  & 7  & -3  & -2  & -3  & -2  & -2  & -3  & $B_1$ \\
\bf 45  & 7  & -3  & -2  & -2  & -3  & -3  & -2  & $B_1$ \\
\bf 46  & 7  & -3  & -2  & -2  & -3  & -2  & -3  & $B_1$ \\
\bf 47  & 7  & -3  & -2  & -2  & -2  & -3  & -3  & $B_1$ \\
\bf 48  & 7  & -2  & -3  & -3  & -3  & -2  & -2  & $B_2$ \\
\bf 49  & 7  & -2  & -3  & -3  & -2  & -3  & -2  & $B_2$ \\
\bf 50  & 7  & -2  & -3  & -3  & -2  & -2  & -3  & $B_2$ \\
\end{tabular}
&
\begin{tabular}{rrrrrrrrr}
\bf \# & $L$ & $E_1$ & $E_2$ & $E_3$ & $E_4$ & $E_5$ & $E_6$ & Divisor Class\\
\hline
\bf 51  & 7  & -2  & -3  & -2  & -3  & -3  & -2  & $B_2$ \\
\bf 52  & 7  & -2  & -3  & -2  & -3  & -2  & -3  & $B_2$ \\
\bf 53  & 7  & -2  & -3  & -2  & -2  & -3  & -3  & $B_2$ \\
\bf 54  & 7  & -2  & -2  & -3  & -3  & -3  & -2  & $B_3$ \\
\bf 55  & 7  & -2  & -2  & -3  & -3  & -2  & -3  & $B_3$ \\
\bf 56  & 7  & -2  & -2  & -3  & -2  & -3  & -3  & $B_3$ \\
\bf 57  & 7  & -2  & -2  & -2  & -3  & -3  & -3  & $B_4$ \\
\bf 58  & 6  & -3  & -2  & -2  & -2  & -2  & -1  & $B_1$ \\
\bf 59  & 6  & -3  & -2  & -2  & -2  & -1  & -2  & $B_1$ \\
\bf 60  & 6  & -3  & -2  & -2  & -1  & -2  & -2  & $B_1$ \\
\bf 61  & 6  & -3  & -2  & -1  & -2  & -2  & -2  & $B_1$ \\
\bf 62  & 6  & -3  & -1  & -2  & -2  & -2  & -2  & $B_1$ \\
\bf 63  & 6  & -2  & -3  & -2  & -2  & -2  & -1  & $B_2$ \\
\bf 64  & 6  & -2  & -3  & -2  & -2  & -1  & -2  & $B_2$ \\
\bf 65  & 6  & -2  & -3  & -2  & -1  & -2  & -2  & $B_2$ \\
\bf 66  & 6  & -2  & -3  & -1  & -2  & -2  & -2  & $B_2$ \\
\bf 67  & 6  & -2  & -2  & -3  & -2  & -2  & -1  & $B_3$ \\
\bf 68  & 6  & -2  & -2  & -3  & -2  & -1  & -2  & $B_3$ \\
\bf 69  & 6  & -2  & -2  & -3  & -1  & -2  & -2  & $B_3$ \\
\bf 70  & 6  & -2  & -2  & -2  & -3  & -2  & -1  & $B_4$ \\
\bf 71  & 6  & -2  & -2  & -2  & -3  & -1  & -2  & $B_4$ \\
\bf 72  & 6  & -2  & -2  & -2  & -2  & -3  & -1  & $B_5$ \\
\bf 73  & 6  & -2  & -2  & -2  & -2  & -1  & -3  & $B_6$ \\
\bf 74  & 6  & -2  & -2  & -2  & -1  & -3  & -2  & $B_5$ \\
\bf 75  & 6  & -2  & -2  & -2  & -1  & -2  & -3  & $B_6$ \\
\bf 76  & 6  & -2  & -2  & -1  & -3  & -2  & -2  & $B_4$ \\
\bf 77  & 6  & -2  & -2  & -1  & -2  & -3  & -2  & $B_5$ \\
\bf 78  & 6  & -2  & -2  & -1  & -2  & -2  & -3  & $B_6$ \\
\bf 79  & 6  & -2  & -1  & -3  & -2  & -2  & -2  & $B_3$ \\
\bf 80  & 6  & -2  & -1  & -2  & -3  & -2  & -2  & $B_4$ \\
\bf 81  & 6  & -2  & -1  & -2  & -2  & -3  & -2  & $B_5$ \\
\bf 82  & 6  & -2  & -1  & -2  & -2  & -2  & -3  & $B_6$ \\
\bf 83  & 6  & -1  & -3  & -2  & -2  & -2  & -2  & $B_2$ \\
\bf 84  & 6  & -1  & -2  & -3  & -2  & -2  & -2  & $B_3$ \\
\bf 85  & 6  & -1  & -2  & -2  & -3  & -2  & -2  & $B_4$ \\
\bf 86  & 6  & -1  & -2  & -2  & -2  & -3  & -2  & $B_5$ \\
\bf 87  & 6  & -1  & -2  & -2  & -2  & -2  & -3  & $B_6$ \\
\bf 88  & 5  & -2  & -2  & -2  & -2  & -2  & -1  & $B_1$ \\
\bf 89  & 5  & -2  & -2  & -2  & -2  & -1  & -2  & $B_1$ \\
\bf 90  & 5  & -2  & -2  & -2  & -1  & -2  & -2  & $B_1$ \\
\bf 91  & 5  & -2  & -2  & -1  & -2  & -2  & -2  & $B_1$ \\
\bf 92  & 5  & -2  & -1  & -2  & -2  & -2  & -2  & $B_1$ \\
\bf 93  & 5  & -1  & -2  & -2  & -2  & -2  & -2  & $B_2$ \\
\bf 94  & 3  & -1  & -1  & -1  & -1  & -1  &  0  & $L_1$ \\
\bf 95  & 3  & -1  & -1  & -1  & -1  &  0  & -1  & $L_1$ \\
\bf 96  & 3  & -1  & -1  & -1  &  0  & -1  & -1  & $L_1$ \\
\bf 97  & 3  & -1  & -1  &  0  & -1  & -1  & -1  & $L_1$ \\
\bf 98  & 3  & -1  &  0  & -1  & -1  & -1  & -1  & $L_1$ \\
\bf 99  & 3  &  0  & -1  & -1  & -1  & -1  & -1  & $L_2$ \\
\\
\end{tabular} \\
\hline
\end{tabular} \\
\end{centering}

\end{document}